\documentclass[onecolumn]{revtex4}

\usepackage{amssymb,amsmath}

\usepackage{natbib}

\begin{document}

\title{The prime analog of the Kepler-Bouwkamp constant}
\author{Adrian R. Kitson}
\email{a.r.kitson@massey.ac.nz}
\affiliation{Institute of Fundamental Sciences, Massey University, Private Bag 11 222, Palmerston North, New Zealand}

\maketitle

Begin with a circle of radius \(R_{1}\) and circumscribe it with an equilateral triangle. Circumscribe the triangle with another circle. Let the radius of the second circle be \(R_{2}\). Circumscribe the second circle with a square. Continue in this manor, circumscribing with a circle, then a regular pentagon, circle, regular hexagon, \emph{ad infinitum}, each time adding one more side to the regular polygon. The limit of the ratio of the radius of the outer circle to the inner is
\begin{displaymath}
K=\lim_{n\rightarrow\infty}\frac{R_{n}}{R_{1}}=8.7000366252\ldots,
\end{displaymath}
which is sometimes referred to as the \emph{polygon circumscribing constant}~\cite{sloane}. Since
\begin{displaymath}
\frac{R_{n}}{R_{1}}=\prod_{k=3}^{n}\sec\left(\frac{\pi}{k}\right),
\end{displaymath}
we have
\begin{displaymath}
K=\prod_{k=3}^{\infty}\sec\left(\frac{\pi}{k}\right).
\end{displaymath}
If, instead of circumscribing, we inscribed at each step we get the \emph{Kepler-Bouwkamp constant}~\cite{finch}, which is
\begin{displaymath}
\rho=\frac{1}{K}=\prod_{k=3}^{\infty}\cos\left(\frac{\pi}{k}\right)=0.1149420448\ldots.
\end{displaymath}
What if only regular polygons with prime (greater than 2) number of edges are considered? In other words, what is
\begin{equation}
\label{eq:kp}
K_{\mbox{p}}=\prod_{p\ge 3}\sec\left(\frac{\pi}{p}\right),
\end{equation}
where \(p\) denotes the primes \(p=2,3,5,\ldots\)?

Taking the natural log of both sides of equation~(\ref{eq:kp}),
\begin{equation}
\label{eq:lnkp}
\ln\left(K_{\mbox{p}}\right)=\sum_{p\ge 3}\ln\left(\sec\left(\frac{\pi}{p}\right)\right).
\end{equation}
Expanding the summand,
\begin{equation}
\label{eq:lnsec}
\ln\left(\sec\left(\frac{\pi}{p}\right)\right)=\sum_{k=1}^{\infty}\frac{2^{2k}\left(2^{2k}-1\right)|B_{2k}|}{2k\left(2k\right)!}\left(\frac{\pi}{p}\right)^{2k},
\end{equation}
where \(B_{n}\) is the \(n\)th Bernoulli number, which can be obtained from
\begin{displaymath}
\frac{x}{\exp\left(x\right)-1}=\sum_{n=0}^{\infty}\frac{B_{n}x^{n}}{n!}.
\end{displaymath}
Using the identity
\begin{displaymath}
B_{2n}=\frac{\left(-1\right)^{n-1}2\left(2n\right)!}{\left(2\pi\right)^{2n}}\,\zeta\left(2n\right),
\end{displaymath}
where \(\zeta\) is the Riemann zeta function, equation~(\ref{eq:lnsec}) becomes
\begin{displaymath}
\ln\left(\sec\left(\frac{\pi}{p}\right)\right)=\sum_{k=1}^{\infty}\frac{\left(2^{2k}-1\right)\zeta\left(2k\right)}{k}\frac{1}{p^{2k}}.
\end{displaymath}
Thus, equation~(\ref{eq:lnkp}) can be written as a double sum,
\begin{equation}
\label{eq:lnkp2}
\ln\left(K_{\mbox{p}}\right)=\sum_{p\ge 3}\sum_{k=1}^{\infty}\frac{\left(2^{2k}-1\right)\zeta\left(2k\right)}{k}\frac{1}{p^{2k}}.
\end{equation}
Equation~(\ref{eq:lnkp2}) is absolutely convergent; the order of summation can be interchanged, \emph{viz}.
\begin{displaymath}
\ln\left(K_{\mbox{p}}\right)=\sum_{k=1}^{\infty}\frac{\left(2^{2k}-1\right)\zeta\left(2k\right)}{k}\sum_{p\ge 3}\frac{1}{p^{2k}}.
\end{displaymath}
The summation over the primes can be written in terms of the prime zeta function, which is
\begin{displaymath}
{\mathcal P}\left(n\right)=\sum_{p}\frac{1}{p^{n}}.
\end{displaymath}
Thus,
\begin{equation}
\label{eq:lnkp4}
\ln\left(K_{\mbox{p}}\right)=\sum_{k=1}^{\infty}\frac{\left(2^{2k}-1\right)\zeta\left(2k\right)}{k}\left({\mathcal P}\left(2k\right)-\frac{1}{2^{2k}}\right).
\end{equation}
Perhaps the only unfamiliar term is the prime zeta function, but this can be evaluated in terms of the more familiar Riemann zeta function using the identity~\cite{titchmarsh}
\begin{displaymath}
{\mathcal P}\left(n\right)=\sum_{k=1}^{\infty}\frac{\mu\left(k\right)\ln\left(\zeta\left(n k\right)\right)}{k},
\end{displaymath}
where \(\mu\) is the M\"{o}bius function.

Equation~(\ref{eq:lnkp4}) quickly converges; taking the exponential gives
\begin{displaymath}
K_{\mbox{p}}=3.1965944300\ldots.
\end{displaymath}
The prime analog of the Kepler-Bouwkamp constant is
\begin{displaymath}
\rho_{\mbox{p}}=\frac{1}{K_{\mbox{p}}}=\prod_{p\ge 3}\cos\left(\frac{\pi}{p}\right)=0.3128329295\ldots.
\end{displaymath}
We also have, with no extra work, the following products over the positive nonprimes, \(q=1,4,6,\ldots\),
\begin{displaymath}
K_{\mbox{q}}=\frac{K}{K_{\mbox{p}}}=\prod_{q\ge 4}\sec\left(\frac{\pi}{q}\right)=2.7216579443\ldots,
\end{displaymath}
and
\begin{displaymath}
\rho_{\mbox{q}}=\frac{1}{K_{\mbox{q}}}=\prod_{q\ge 4}\cos\left(\frac{\pi}{q}\right)=0.3674231004\ldots.
\end{displaymath}

\end{document}